\documentclass{conm-p-l}
\usepackage{latexsym}
\usepackage{amsmath}
\usepackage{amsthm}
\usepackage{amsfonts}
\usepackage{amscd}
\usepackage{amsxtra}
\title{Hypercyclic Differentiation Operators}

\author{Richard Aron}
\address{Department of Mathematics and Computer Science \\
         Kent State University \\
         Kent, OH 44242. U.S.A.}
\email{aron@mcs.kent.edu}

\author{Juan B\`es}
\address{Department of Mathematics and Statistics \\
         Bowling Green State University \\
         Bowling Green, OH 43403. U.S.A.}
\email{jbes@bgnet.bgsu.edu}

\date{August 1998}

\begin{document}
\begin{abstract}
A classical theorem due to G. D. Birkhoff states that there
exists an entire function whose translates approximate
any given entire function, as accurately as desired, over
any ball of the complex plane. We show this result may be
generalized to the space $H_{bc}(E)$ of entire functions
of compact bounded type defined on a Banach space $E$ with
separable dual.

\end{abstract}
\maketitle

\newcommand{\K}{\rm I\kern-.19emK}
\newcommand{\C}{\mathbb{C}}
\newcommand{\R}{\mathbb{R}}
\newcommand{\N}{\mathbb{N}}\newcommand{\Z}{\mathbb{Z}}

\newcommand{\e}{\epsilon }
\newcommand{\f}{\varphi  }
\newcommand{\ti}{\tilde  }
\newcommand{\sn}{\sum_{n=0}^{\infty }}


\newtheorem{theorem}{Theorem}[section]
\newtheorem{lemma}[theorem]{Lemma}
\newtheorem{corollary}[theorem]{Corollary} 

\theoremstyle{definition}
\newtheorem{definition}[theorem]{Definition}
\newtheorem{example}[theorem]{Example}

\theoremstyle{remark}
\newtheorem{remark}[theorem]{Remark}

\numberwithin{equation}{section}


%
%

\section{{\bf Introduction.}}
The Invariant Subspace Problem asks whether every bounded linear operator on an infinite
dimensional Banach space admits a non-trivial, closed invariant subspace, that is, a closed linear
manifold which is neither $\{ 0 \}$ nor the whole space and which is mapped into itself
by the operator. A related question was raised by M. Edelstein and H. Radjavi:
Does every such operator have a non-trivial closed invariant {\em set}? (\cite{Ki82}, p. 1)
\medskip

P. Enflo created an operator on a (non-Hilbertian) Banach space without non-trivial closed
subspaces, to show the answer is negative in general \cite{En87}. 
Later, C. Read constructed an operator on $l_1$ with the same property, and also an operator
on a Banach space admitting no non-trivial closed invariant subsets (\cite{Re85}, \cite{Re88}).
The Invariant Subspace Problem remains open for the case in which the Banach space is a Hilbert space.
\medskip

The notions of cyclic and hypercyclic vectors arise naturally in this study of invariant
subspaces and subsets.
Suppose $T$ is a continuous linear operator on an $\mathcal{F}$-space (i.e., a complete, linear
metric space) $X$. A vector $x_0\in X$ is {\bf cyclic} for $T$ provided the linear span of the
orbit
\[
\mbox{Orb}\{x_0,T\}=\{ x_0, Tx_0, T^2x_0,\dots \}
\]
is dense in $X$. If the orbit itself is dense, the operator $T$ and associated vector
$x_0$ are called {\bf hypercyclic}. In this way, an operator $T$ will lack non-trivial, closed invariant subspaces
(subsets) if and only if {\em all} the non-zero vectors are cyclic (hypercyclic) for~$T$.
\medskip

One of the roots of the modern study of hypercyclicity comes from an intriguing
observation of G.D. Birkhoff concerning the orbits of translation operators 
acting on the space of entire functions \cite{Bi29}. Let $H(\C )$ denote 
the space of entire functions 
of one complex variable, endowed with the topology of uniform convergence on 
compact subsets
of the plane, and let $\tau_b :H(\C )\rightarrow H(\C )$ be the operator 
of ``translation by the complex number $b$'' defined by
\[
\tau_b(f)(z):=f(z+b) \ \ \ \ \ \ \ \ (f\in H(\C ), z\in \C )
\]
Birkhoff's theorem asserts that if $b\ne 0$, then the operator $\tau_b$ is hypercyclic.
\medskip

In other words, Birkhoff's theorem gives a ``universal'' entire function 
which, over any compact set, has translates that approximate 
any entire function as accurately as desired. In the language of dynamics, it 
states that each translation operator is {\em topologically transitive} on $H(\C )$, and this
is a major step in proving that these operators are actually {\em chaotic} in one of the commonly
accepted senses (cf. \cite{GoS91}, section 6).
\medskip

Birkhoff's operator may be viewed in several different contexts. 
 As a translation 
operator, W. Seidel and J. Walsh provided a non-Euclidean version of Birkhoff's theorem 
in the space of holomorphic functions on the unit disc, by replacing translation by 
certain conformal disc automorphisms \cite{SeW41}.
\smallskip

It is essentially the only hypercyclic composition operator on $H(\C )$.
Given  $\f \in H(\C )$ fixed, the operator $f  \mapsto f \circ \f $ on $H(\C )$
 is hypercyclic if and only if $\f (z) = z+b$ for some
non-zero complex number $b$ (\cite{jptesis}, Corollary 2.3). 
However, in other spaces such as the Hardy space $H^2$ the situation is quite
different. P. Bourdon and J. Shapiro have studied the cyclic and hypercyclic
behaviour of composition operators in this space, obtaining complete results
for linear fractional maps of the unit disc into itself, and using these maps
as ``models'' to extend this classification to a large class of
composition operators~(\cite{BoS97}, \cite{Sh93}).
\smallskip

Birkhoff's translation has also been regarded as a differentiation 
operator. That is, $\tau_b=e^{bD}$,
so that for every $f\in H(\C )$,
\[
\tau_b (f) = \sum_{n=o}^{\infty} \frac{b^n}{n!}\, D^nf  \ \ \ \ \ \ \ \
\ \left( \mbox{{\rm convergence in $H(\C )$}} \right) .
\]
From this perspective, G. R. MacLane showed in 1952 that the differentiation
operator
\[ 
 f\mapsto Df
\]
 is also hypercyclic on $H(\C )$. Thus, MacLane's theorem assures that 
there exists an entire function
whose sequence of derivatives has {\em every } entire function as a limit point
(\cite{Ma52}, Theorem 6). 

In 1991, a beautiful generalization of Birkhoff's and MacLane's theorems was
provided by G. Godefroy and J. Shapiro (\cite{GoS91}, Thm 5.1):
\begin{theorem} {\bf (Godefroy-Shapiro)} \label{T:GoS91}
Let $\Phi (z)=\sum_{|\alpha |\ge 0 } c_{\alpha} z^{\alpha}$ be a non-constant,
entire function on $\C^N$of exponential type. Then the operator
\[
f\overset{\Phi (D)}{\longmapsto }\sum_{|\alpha |\ge 0 } c_{\alpha} D^{\alpha}f \ \ \ \ f\in H(\C^N )
\]
is hypercyclic.
\end{theorem} 

Moreover, they showed that a continuous linear operator $T$ on $H(\C^N)$ commutes with
translations if and only if it commutes with the differentiation operators
$\frac{\partial}{\partial z_k} \ (1\le k\le N)$, and if and only if it is of the form
$T=\Phi (D)$, for some $\Phi\in H(\C^N )$ of exponential type (\cite{GoS91}, Prop. 5.2).
\smallskip

 In the next section we will consider Theorem~\ref{T:GoS91} for the case $N=1$, 
and show how it may be 
generalized by replacing $H(\C )$ by the
space $H_{bc}(E)$ of complex valued, entire functions on $E$ of compact 
type that are bounded on bounded subsets of $E$. In particular, an extension of
Birkhoff's theorem to all such spaces will be obtained.
\smallskip

Our proof involves a technique that has become standard in this area, the so-called
Hypercyclicity Criterion. In 1982, C. Kitai isolated conditions that ensure an operator
to be hypercyclic (\cite{Ki82}, Thm 1.4). 
This result was never published, and
it was later rediscovered in a broader form by R. Gethner and J. Shapiro 
(\cite{GeS87},Thm 2.2), who used it to unify the previously mentioned
theorems of Birkhoff, 
MacLane, and Seidel and Walsh, among others.
It may be stated as follows:
\smallskip

\begin{theorem} {\bf ( Kitai - Gethner - Shapiro )} \label{T:KGS} 
Let $X$ be an $\mathcal{F}$-space, and $T:X\rightarrow X$ be linear, continuous.
Suppose there exist $X_0, Y_0$ dense subsets of $X$, a sequence $(n_k)$ of positive
integers, and a sequence of mappings
(possibly nonlinear, possibly not continuous) $S_{n_k}: Y_0 \rightarrow X$
so that
\begin{equation}
\begin{align}
i)\phantom{ii}\ & T^{n_k}
\underset {k\to \infty }{\rightarrow 0} \ \mbox{pointwise on } X_0. \ \ \ \ \ \ \ \ \ \ \notag \\
ii)\phantom{i}\ & S_{n_k}
\underset {k\to \infty }{\rightarrow 0} \ \mbox{pointwise on} \ Y_0. \ \ \ \ \ \ \  \ \ \ \notag \\
iii)\ &  T^{n_k}S_{n_k} = \mbox{Identity} \ \ \mbox{ on }Y_0. \ \ \ \ \  \ \  \ \ \ \notag
\end{align}
\end{equation}
Then $T$ is hypercyclic.
\end{theorem}

{\bf Note.} It is an open problem whether the above conditions are also necessary
(\cite{LeM97}, p 544):
\begin{equation} \label{eq:preg}
\mbox{Does every hypercyclic operator satisfy the Hypercyclicity Criterion? }
\end{equation}

There are hypercyclic operators, though, that do not satisfy the Criterion for the
{\em entire } sequence $(n_k)=(k)$ of positive integers. Solving a problem of D. Herrero, H. Salas constructed a weighted shift $A$ so that
both it and its Hilbert Adjoint $A^*$ were hypercyclic, while their direct sum $A\oplus A^*$
was not. As a 
consequence, $A$ and $A^*$ 
could not simultaneously satisfy the Criterion for the
sequence $(n_k)=(k).$  For a while, some people were 
misled into thinking that it was a
counterexample to (\ref{eq:preg}) (cf. \cite{He91}, p97), but both $A$ and $A^*$ were later noticed in (\cite{LeM98}, Section 2) to satisfy the Criterion (cf. also \cite{jptesis}, Note 1.22). 

If requiring only pointwise convergence in (iii) instead of equality,  it is shown in \cite{BeP98}
that
question (\ref{eq:preg}) is equivalent to the following one of D. Herrero (\cite{He91}, p.~97):
\[
\mbox{Is $T\oplus T$ hypercyclic whenever $T$ is?}
\]
(A direct sum $A\oplus B$ may fail to be hypercyclic even if each of
$A$ and $B$ is (\cite{Sa95}, Corollary 2.6)).

Interest in answering these questions comes also from the study of hypercyclic subspaces,
that is, linear subspaces where all the non-zero vectors are hypercyclic. Although
every hypercyclic operator admits a dense, invariant hypercyclic subspace
(\cite{He91}, \cite{Bo93}, \cite{Bes98}, \cite{Be87}), 
they do not always admit a {\em closed} and infinite dimensional 
one. The situation is well understood for operators satisfying the Criterion.
For instance, an operator on a Banach space satisfying the Hypercyclicity Criterion
will admit a closed, infinite dimensional hypercyclic subspace if and only if
its essential spectrum intersects the closed unit disk (\cite{LeM98}, Thm 2.1; cf also
\cite{Cha98}, \cite{Mo96}, \cite{LeM97}, \cite{BeM95}).

%
%
%
%

%
%

\section[Birkhoff's theorem for the space $H_{bc}(E)$.]{{\bf Birkhoff's theorem 
for the space $H_{bc}(E)$. }}
Given a Banach Space E, let $H_{bc}(E)$ \ be the 
Fr\'echet Algebra generated
 by the elements of
the dual space
 $E^{\ast }$, endowed with the topology of uniform
 convergence on balls of E.
The space $H_{bc}(E)$ \ consists of entire functions $ f:\, E \rightarrow \C$
of so-called compact type that are bounded on bounded subsets of $E$.
That is, $H_{bc}(E)$ consists of all functions of the form
\[ `
f=\sn P_n ,
\]
where
\begin{equation} \notag
\cases \  P_n \in
\overline{\text{span}} \left \{ \f^n : \ \f \in E^{\ast } \right \} \ \ \ n=0,1,2,\dots \\
\| P_n \|^{\frac1n} =
\left(\,\text{sup}_{\|x\|\le 1} |P_n(x)|\, \right)^{\frac1n} 
\underset{n\rightarrow \infty }{\rightarrow } 0.
\endcases
\end{equation}

Parenthetically, we remark that the terminology {\em compact
type} apparently comes from the fact that if $E^{\ast}$ has the
approximation property (see, e.g., 
\cite{M}, p. 195), then any compact linear mapping
$E \to E^{\ast}$ is a uniform limit of finite rank mappings,
which are of the form $\sum_{j = 1}^k \phi_j \psi_j.$ The associated
$2-$homogeneous polynomial, $x \in E \to 
\sum_{j = 1}^k \phi_j(x)\psi_j(x)$ is of {\em finite type}, and
limits of such finite type polynomials are called $2-$homogeneous
polynomials of compact type. Of course, the situation generalizes
to $n-$homogeneous polynomials. Also, the
space $H_{bc}(E)$ is not as arcane as may first appear to the
non-specialist. Indeed, if $E^{\ast}$ has the approximation property,
then $H_{bc}(E)$ has a useful alternative description, as
the space of entire functions which are uniformly continuous on bounded
subsets of $E,$ with respect to the weak topology (\cite{A}, \cite{Di}). 
Moreover if the Michael problem, which asks whether 
every $\C-$valued homomorphism on a complex Fr\'echet algebra 
is automatically continuous, has an affirmative
solution for the Fr\'echet algebra $H_{bc}(c_0),$ then it has an
affirmative solution for every Fr\'echet algebra (\cite{M}).

\smallskip

Given $a\in E$ fixed and $\Phi = \sn c_n z^n \,\in H(\C )$ of exponential
type, let $\Phi_a(D) : H_{bc}(E) \rightarrow H_{bc}(E)$ be defined as
\begin{equation}
\Phi_a(D) (f) = \sn c_n \widehat{d^n f(.)}a\, . \notag 
\end{equation}
That is,
\begin{equation} \notag 
\Phi_a(D) (f) (x) = \sn c_n \widehat{d^nf(x)}(a) \ \ \ \text{for all }
\ x\in E\, ,
\end{equation}
where $\widehat{d^nf(x)}$ is the n-homogeneous polynomial
associated with the $\text{n}^{th}$-Fr\'echet derivative of $f$ at the point 
$x\in E$. The operator $\Phi_a(D)$ is continuous on $H_{bc}(E)$, and when
$\Phi (z) = e^z$ the induced operator is 
``translation by $a$'' (\cite{Gu66}, Lemma 8). i.e.,
\begin{equation} \label{eq:Ta}
\tau_a(f)(x)= f(x+a)= \sn \frac1{n!} \widehat{d^n f(x)}(a) = \Phi_a(D) (f) (x)
\end{equation}
for all $x\in E$ \ and all $ f\in H_{bc}(E)$.
We are ready now to state the main result:

\begin{theorem} \label{T:ABa}
Let E be a Banach space with separable dual $E^{\ast }$, and $0\ne a\in E$.
Let $\Phi (z) = \sn c_n z^n \in H(\C )$ be non-constant, of exponential type.
Then the operator 
\begin{align}
\Phi_a(D): H_{bc} &(E) \longrightarrow H_{bc}(E)   \notag \\
                f    &\mapsto  \sn c_n \widehat{d^nf(.)}a \notag
\end{align}
is hypercyclic.
\end{theorem}

{\bf Note.} If $E^*$ fails to be separable, so does $H_{bc}(E)$, and this prevents
$\Phi_a(D)$ from being hypercyclic. Thus, from now on we assume $E^*$ to be separable.

In particular, Birkhoff's theorem for the spaces $H_{bc}(E)$ follows:

\begin{corollary} 
Let E be a Banach space with separable dual $E^{\ast }$, and $0\ne a\in E$.
Then the operator ``translation by `a' '' 
\[
\tau_a(f)(x)= f(x+a) \ \ \ \ \ \ (x\in E)
\]
is hypercyclic on $H_{bc}(E)$.
\end{corollary}

\begin{proof} It follows immediately from 
   (\ref{eq:Ta}) and Theorem~\ref{T:ABa}.
\end{proof}

We will make use of the following two lemmas. 
\begin{lemma} \label{L:lemma1} \mbox{}
$
\mathcal{B} = \left \{ \ e^{\varphi } : \ \varphi \in E^* \right \}
$ \ is a linearly independent subset of $H_{bc}(E)$.
\end{lemma}

\begin{proof}
Let $\left \{ e^{\f_i} \right \}_{i\in I}$ be a maximal linearly
independent subset of $\mathcal B$.
Fix $\f \in E^{\ast }$, and assume there exist non-zero constants
$c_{i_1},\dots , c_{i_r}\in \C$ so that
\begin{equation}
c_{i_1}e^{\f_{i_1}}+\dots +c_{i_r}e^{\f_{i_r}} = e^{\f } \label{eq:dos2}
\end{equation}

Let $a\in E$ be arbitrary. Applying  the operator
$
f \mapsto df(.)a 
$ in (\ref{eq:dos2}), it follows that
\begin{equation}
c_{i_1}\f_{i_1} (a)e^{\f_{i_1}}+\dots +c_{i_r}\f_{i_r} (a)e^{\f_{i_r}}
= \f (a)e^{\f }. \label{eq:aaa}
\end{equation}

Since $\left \{ e^{\f_i} \right \}_{i\in I}$ is linearly independent and
$c_{i_1},\dots , c_{i_r}$ are non-zero, by (\ref{eq:dos2}) and
(\ref{eq:aaa}) we have
\[\f_{i_1} (a) =\dots = \f_{i_r} (a) = \f (a).\]
Since $a\in E$ is arbitrary, \[ \f_{i_1} =\dots = \f_{i_r} = \f .\]
Hence, $\left \{ \f_i \right \}_{i\in I} = E^{\ast }$,
and Lemma \ref{L:lemma1} follows.
\end{proof}

	Our next lemma generalizes a result of Gupta,
which was needed in \cite{Gu66} to obtain information about
the range and kernel of convolution operators.
\begin{lemma} \label{L:lemma4} \mbox{}
Let $U$ be a non-empty open subset of $E^*$. Then
\begin{equation} \notag
S=\text{span} \left \{ e^{\f} : \ \f \in U \right \}
\end{equation}
is dense in $H_{bc}(E)$.
\end{lemma}
\begin{proof} 
\ Given $\f_0 \in E^*$, the map
\begin{equation} \notag
\underset{g \mapsto e^{\f_0}g}{H_{bc}(E) \longrightarrow H_{bc}(E) }
\end{equation}
is a homeomorphism. So
\[ \overline{\mbox{span}}\left \{ e^{\f_0 +\f} : \ \f \in U \right \}=H_{bc}(E)
\Leftrightarrow \overline{\mbox{span}}\left \{ e^{\f} : \ \f \in U \right \} =H_{bc}(E).
\]
Hence, we may assume that $0\in U$. Reducing $U$ if necessary,
we may also assume that for some $\delta >0$,
\begin{equation} \notag
U = \left \{ \f \in E^* \, : \ \| \f \| < \delta \right \}.
\end{equation}
In particular, $1\in \overline S$,\,  for $0\in U$. It will then suffice to show
the following:
\smallskip

{\em Claim: \ \  For every $\f \in U$ and $ n\ge 1$, \ 
$\f^n \in \overline S $}.
\smallskip

To prove the claim, let $\f \in U$ and suppose the claim is true for $n\le k-1$.
Then for each $0<t<1$ we have 
\begin{equation} \notag
g_t= \frac{e^{t\f} -1 -t\f - \frac{[t\f]^2}{2!} - \dots - 
\frac{[t\f]^{k-1}}{(k-1)!}}{t^k} \in
\overline S,
\end{equation}
since $t\f \in U$. So given $x\in E$,

\begin{align}
\left | { \left ( g_t-\frac{\f^k}{k!} \right ) (x) } \right | &=
\left | \, { 
\frac1{t^k} \left [ e^{t\f} -1 -t\f - \frac{[t\f]^2}{2!} - \dots - \frac{[t\f]^k}{k!} \right ] (x) }\, \right |
\notag  \\
&\le  \left |\, \frac1{t^k} \sum_{n\ge k+1} \frac{[t\f]^n}{n!}(x)\, \right | \notag \\
&\le \, t \, \sum_{n\ge k+1} t^{n-k-1} \frac{|\f (x)|^n}{n!} \notag \\
&\le \, t \ e^{\delta \|x \|}. \notag
\end{align}

Thus,  $ g_t \underset{t \to 0}{\rightarrow} \frac{\f^k}{k!} $ 
\ in $H_{bc}(E)$, and $\frac{\f^k}{k!} \in \overline S$. So the claim holds. 
\end{proof}

\begin{proof}[Proof of Theorem \ref{T:ABc}]
Consider the function $g:E^* \longrightarrow \C $ defined by 
\[
g(\f )=\sum_{n=0}^{\infty } P_n (\f ) ,
\]
where each $P_n:E^* \longrightarrow \C$ is the $n$-homogeneous polynomial defined by
\[
P_n (\f ) = c_n \f^n (a) \ \ \ (n \ge 0).
\] 
 Now, since $\Phi $ is of exponential type, there exists $R>0$ so that
\[
|c_n | \le \frac{R^n}{n!}
\ \ \ \ (n\ge 1).
\]
Given $\f \in E^*$ with $\| \f \| \le 1 $ and $n\ge 1$, 
\[
|P_n (\f ) | = |c_n | \, |\f^n (a) | \le \frac{R^n}{n!} \| \f \|^n \|a \|^n 
\le \left ( \frac{ e\, \|a\| R }{n}  \right )^n, 
\]
and so
\begin{equation} \notag
\| P_n \|^{\frac1n } \underset{n\to \infty }{\longrightarrow} 0.
\end{equation}
  
Thus, $g:E^* \longrightarrow \C$ is entire (moreover, it is of bounded type),
and non-constant, since $\Phi $ is non-constant.
So the sets
\begin{align}
U& := \left \{ \, \f \in E^* : \  |\sum_{n=0}^{\infty }
c_n \f^n (a) | = |g(\f )| < 1 \ \right \} \notag \\
V &:=
\left \{ \, \f \in E^* : \  |\sum_{n=0}^{\infty } c_n \f^n (a) | = |g(\f )| > 1 \ \right \}
 \notag
\end{align}
are both open, non-empty. Hence, according to Lemma \ref{L:lemma4},
\begin{equation} 
\begin{align}
X_0 &:= \text{span}\left \{ e^{\f} : \  \f\in U \, \right \}  \label{seq:veintidos1} \\
Y_0 &:= \text{span}\left \{ e^{\f} :  \ \f\in V \, \right \}  \label{seq:veintidos2} 
\end{align}
\end{equation}
are both dense subspaces of $H_{bc}(E)$.
Next, notice that if $T=\Phi_a(D)$, 
given $\f\in E^*$ 
\begin{align}
T(e^{\f }) &= \sum_{n=0}^{\infty } c_n \widehat{d^n(e^{\f} )}a  \notag \\
&=\sum_{n=0}^{\infty } c_n \f^n (a) \ e^{\f } \notag \\
&= g(\f ) e^{\f }. \notag
\end{align}
By (\ref{seq:veintidos1}), 
\begin{equation} \notag
T^n\underset{n\to \infty }{\longrightarrow} 0 \ \ \text{pointwise on } \ X_0.
\end{equation}
Also, by Lemma \ref{L:lemma1} there exists a (possibly discontinuous) linear map 
$S:Y_0 \longrightarrow Y_0$ determined by
\begin{equation} \label{eq:veinticuatro24}
S(e^{\f }) = [ g(\f )]^{-1} \ e^{\f } \ \ \ \ \ (\f\in E^*)
\end{equation}
which by (\ref{seq:veintidos2}) and (\ref{eq:veinticuatro24}) satisfies 
\[
\begin{cases}
 S^n \underset{n\to \infty }{\longrightarrow} 0 \ \ \text{pointwise on } \ Y_0.\\
 TS = id_{Y_0}  \ \ \ \ \ \ \text{ on } Y_0.
\end{cases}
\] 
By Theorem \ref{T:KGS}, $T=\Phi_a(D)$ is hypercyclic. 
\end{proof}

{\bf Remarks:}

i) For $E=\C$ and $a=1$, Theorem~\ref{T:ABa} yields Theorem~\ref{T:GoS91}
(for the case $N=1$), since $H_{bc}(\C )=H(\C )$, and given any $\Phi\in H(\C )$ of
exponential type, $\Phi (D)=\Phi_1(D)$.
\medskip

ii) Using similar arguments, one may obtain the following two theorems (cf \cite{jptesis}, Chapter 2):

\begin{theorem} \label{T:ABb}
Let E be a Banach space with separable dual $E^{\ast }$.
Let $\{ b_1, b_2, \dots , b_r \}$ be a linearly independent subset of $E$,
and let $\Phi^1,\dots ,\Phi^r$ in $H(\C )$ be non-constant, of
exponential type.
Then
\begin{equation} \notag
A=\Phi^1_{b_1}(D) + \Phi^2_{b_2}(D) + \dots + \Phi^r_{b_r}(D)
\end{equation}
is hypercyclic on $H_{bc}(E)$.
\end{theorem}

\begin{theorem} \label{T:ABc}
Let E be a Banach space with separable dual $E^{\ast }$, and let
 $\{ b_i \}_{i\ge 0}\subset E\setminus \{ 0\}$ be bounded.
Then the operator
\begin{align}
A: H_{bc}(&E) \longrightarrow H_{bc}(E)  \notag  \\
     f&\mapsto \sn c_n  \widehat{d^nf(.)}b_n \notag
\end{align}
is hypercyclic, where $\Phi (z) = \sn c_n z^n \in H(\C )$ is entire, non-constant
and of exponential type.
\end{theorem}
\medskip


%
%


\begin{thebibliography}{99}

\bibitem{A}
R. M. Aron, {\em Weakly uniformly and weakly sequentially continuous
entire functions,} in {\em Adv. in Holomorphy} (J. A. Barroso, ed.),
North-Holland (1979), 47 - 65.

\bibitem{Be87}
B. Beauzamy,
{\em An operator in a separable Hilbert space with many hypercyclic vectors,}
Studia Math.\ {\bf 87} (1987), 71-78.

\bibitem{BeM95}
L. Bernal-Gonz\'alez and A. Montes-Rodr\'\i guez,
{\em Non-finite dimensional closed vector spaces of universal functions for 
composition operators }
J. Approx. Theory\ {\bf 82}, (1995), 375-391.

\bibitem{Bes98}
J. B\`es,
{\em Invariant manifolds of hypercyclic vectors for the real scalar case, }
Proc. Amer. Math. Soc.\ {\bf 127} (1999) pp1801-1804.

\bibitem{jptesis}
J. B\`es, {\em Three problems on hypercyclic operators,} Thesis, Kent State Univ.,
(1998).

\bibitem{BeP98}
J. B\`es and A. Peris-Manguillot,
{\em Hereditarily hypercyclic operators,} preprint.

\bibitem{Bi29}
G. D. Birkhoff,
{\em D\'emonstration d'un th\'eoreme \'el\'ementaire sur les
fonctions enti\`eres,}
 C. R. Acad. Sci. Paris\ {\bf 189} (1929), 473-475.

\bibitem{Bo93}
P. Bourdon,
{\it Invariant Manifolds of Hypercyclic Vectors},
Proc. Amer. Math. Soc. \ {\bf 118} No. 3 (1993), 845-847.

\bibitem{BoS97}
P. S. Bourdon and J. H. Shapiro, 
{\em Cyclic phenomena for composition operators,}
Mem. Amer. Math. Soc.\ {\bf 125} (1997).
 
\bibitem{Cha98}
K. C. Chan,
{\em Hypercyclicity of the Operator Algebra for a separable Hilbert space, }
preprint.

\bibitem{De89}
R. L. Devaney,
{\em ``An introduction to Chaotic Dynamical Systems''},
$2^{nd}$ ed., Addison-Wesley, Reading, MA. 1989.

\bibitem{Di}
S. Dineen, {\em  Entire functions on $c_{0},$} 
J. Funct. Anal. {\bf 52} no. 2, (1983) 205--218.

\bibitem{En87} P.\ Enflo, 
{\em On the invariant subspace problem for Banach spaces,}
Acta Math.\ {\bf 158} (1987), 213-313.

\bibitem{GeS87} R.\ M.\ Gethner and J.\ H.\ Shapiro, 
{\em Universal vectors for operators on spaces of holomorphic functions,} 
Proc.\ Amer.\ Math.\ Soc.\ {\bf 100}, No.\ 2 (1987), 281-288.

\bibitem{GoS91} G.\ Godefroy and J.\ H.\ Shapiro, 
{\em Operators with dense, invariant, cyclic vector manifolds,} 
J.\ Funct.\ Anal.\ {\bf 98} (1991), 229-269.

\bibitem{Gu66} C. P. Gupta,
{\em Malgrange Theorem for Nuclearly Entire Functions of Bounded Type on
a Banach Space},
Thesis, University of Rochester (1966).

\bibitem{He91}
D. A. Herrero,
{\em  Hypercyclic Operators and Chaos,}
Journal of Op. Theory\ {\bf 28} (1992), 93-103.

\bibitem{Ki82} C.\ Kitai, 
{\em Invariant Closed Sets for Linear Operators,}
Ph.D.\ thesis, Univ.\ of Toronto, 1982.

\bibitem{LeM97} 
F. Le\'on-Saavedra and A. Montes-Rodr\'\i guez,
{\em Linear structure of hypercyclic vectors}
 J. Funct. Anal.\ {\bf 148}\ No. 2 (1997), 524-545.

\bibitem{LeM98}
F. Le\'on-Saavedra and A. Montes-Rodr\'\i guez,
{\em Spectral Theory and Hypercyclic Subspaces,}
Trans. Amer. Math. Soc.\ (to appear).

\bibitem{Ma52}
G. R. MacLane,
{\em Sequences of derivatives and normal families,}
 J. Analyse Math.\ {\bf 2} (1952), 72-87.

\bibitem{Mo96}
A. Montes-Rodr\'\i guez,
{\em Banach spaces of hypercyclic vectors,}
Michigan Math. J.\ {\bf 43} No.~3 (1996), 419-436.

\bibitem{M}
J. Mujica, {\em ``Complex analysis in Banach spaces,''} Math.
Studies, North-Holland, {\bf 120} (1986).

\bibitem{Re85}
C.\ Read, 
{\em A solution to the invariant subspace problem on the space $l_1$,} 
Bull.\ London Math.\ Soc.\ {\bf 17} (1985), 305-317.  

\bibitem{Re88}
C.\ Read,
{\em The invariant subspace problem for a Class of Banach spaces 2: Hypercyclic Operators},
Israel J. Math. \ {\bf 63}, No.~1, 1998.




\bibitem{Sa95} 
H.\ Salas, 
{\em Hypercyclic weighted shifts,} 
Trans.\ Amer.\ Math.\ Soc.\ {\bf 347}, No.\ 3 (1995), 993-1004.

\bibitem{SeW41} 
W. P. Seidel and J. L. Walsh,
{\em On Approximation by Euclidean and non-Euclidean translates of an Analytic Function,}
Bull. Amer. Math. Soc.\ {\bf 47} (1941), 916-920.
 
\bibitem{Sh93}
J. H. Shapiro,
{\em ``Composition Operators and classical function theory'',}
Springer-Verlag\ (1993).
 
\end{thebibliography}
\end{document}